\documentclass{article}

\usepackage{amsmath}
\usepackage{amssymb}
\usepackage{enumitem}
\usepackage{graphicx}
\usepackage{mathdots}
\usepackage{color}

\def\B{{\boldsymbol B}}
\def\Be{{\boldsymbol B}^{\rm even}}
\def\Bo{{\boldsymbol B}^{\rm odd}}

\def\nice{\displaystyle}

\def\End{{\rm End}}

\def\oM{\overline{\mathcal{M}}}
\def\cM{{\mathcal{M}}}

\def\Z{\mathbb{Z}}
\def\C{\mathbb{C}}
\def\Q{\mathbb{Q}}

\def\qed{{\hfill $\Diamond$}}
\def\b1{{\bf 1}}

\def\Aut{{\rm Aut}}

\def\E{\mathrm{E}}
\def\n{\mathrm{n}}
\def\L{\mathrm{L}}
\def\V{\mathrm{V}}
\def\H{\mathrm{H}}
\def\g{\mathrm{g}}
\def\rarr{\rightarrow}
\def\W{\mathsf{W}}
\def\D{\mathsf{D}}

\def\com{\mathbb{C}}

\def\HH{{\mathcal{H}}}
\def\slc{{\mathsf{sl}}_2(\mathbb{C})}
\def\sl{{\mathsf{sl}}}
\def\ch{{\mathsf{ch}}}

\usepackage{theorem}
\newtheorem{definition}{Definition}
\newtheorem{theorem}[definition]{Theorem}
\newtheorem{example}[definition]{Example}
\newtheorem{proposition}[definition]{Proposition}
\newtheorem{corollary}[definition]{Corollary}
\newtheorem{lemma}[definition]{Lemma}

\newtheorem{question}[definition]{Question}

\newcommand{\Hilb}{\mathsf{Hilb}^m(\mathbb{C}^2)}
\newcommand{\Sym}{\mathsf{Sym}^m}

\newcommand{\CC}{\mathbb{C}}
\newcommand{\cI}{\mathcal{I}}
\newcommand{\T}{\mathsf{T}}
\newcommand{\cF}{\mathcal{F}}
\newcommand{\vac}{v_\emptyset}
\newcommand{\lv}{\left |}

\newcommand{\ZZ}{\mathbb{Z}}

\newcommand{\lang}{\left\langle}
\newcommand{\rang}{\right\rangle}

\newcommand{\zz}{{\mathfrak{z}}}

\newcommand{\MM}{\mathsf{M}}

\title{Cohomological field theory calculations}
\author{Rahul Pandharipande}
\date{February 2018}

\begin{document}

\maketitle

\vspace{-20pt}

\begin{abstract} 
Cohomological field theories (CohFTs) were defined in the mid 1990s
by Kontsevich and Manin to capture the formal
properties of the virtual fundamental class in Gromov-Witten theory.
A beautiful classification result for semisimple CohFTs 
 (via the action of the Givental group) was proven by Teleman in 2012.
The Givental-Teleman classification can be used to explicitly
calculate the full CohFT in many interesting cases not
approachable by earlier methods. 

My goal here is to present an introduction to these ideas
together with a survey of the calculations of the CohFTs obtained from 
\begin{enumerate}
\item[$\bullet$]
Witten's classes on the moduli spaces of $r$-spin curves, 
\item[$\bullet$] Chern characters of the
Verlinde bundles on the moduli of curves,
\item[$\bullet$]  Gromov-Witten classes 
of Hilbert schemes of points of $\com^2$. 
\end{enumerate}
The subject is full
of basic open questions.
\end{abstract}

\setcounter{tocdepth}{1} 
%\tableofcontents

%%%%%%%%%%%%%%%%%%%%%%%%%%%%%%%%%%%%%%%%%%%%%%%%%%%%%%%%%%%%%%%%%%%%%%%%%%%

\setcounter{section}{-1}
\section{Introduction}

\subsection{Moduli of curves} \label{xcxc}
The moduli space $\cM_g$ of complete, nonsingular, irreducible, 
 algebraic curves
over $\com$ of genus $g$ has been a central object in
mathematics since Riemann's work  in the middle of the 19th century.
The Deligne-Mumford compactification
$$\cM_g \subset {\oM}_g$$
by nodal curves was defined almost 50 years ago \cite{DM}.

We will be concerned here with the moduli space of curves with marked points,
$$\cM_{g,n} \subset {\oM}_{g,n}\, , $$
in the stable range $2g-2+n>0$. 
As a Deligne-Mumford stack (or orbifold), ${\oM}_{g,n}$
is nonsingular, irreducible, and of (complex) dimension $3g-3+n$.
There are natural forgetful morphisms
$$p:\oM_{g,n+1} \rightarrow \oM_{g,n}$$
dropping the last marking.

The {\em boundary}{\footnote{Since ${\oM}_{g,n}$ is a
closed nonsingular orbifold, the boundary here is {\em not} in
the sense of orbifold with boundary. 
If $g=0$,
there is no boundary map $q$.}} 
of the Deligne-Mumford compactification is the closed locus parameterizing
curves with a least one node,
$$ \partial{\oM}_{g,n}= {\oM}_{g,n} \setminus \cM_{g,n}\, .$$
By identifying the last two markings of a single $(n+2)$-pointed curve of genus $g-1$,
we obtain a morphism
$$
q : \oM_{g-1, n+2} \to \oM_{g,n}\ .
$$
Similarly, by identifying the last markings of separate pointed curves,
we obtain
$$
r: \oM_{g_1, n_1+1} \times \oM_{g_2, n_2+1} \to \oM_{g,n}\, , 
$$
where $n=n_1+n_2$ and $g=g_1+g_2$. The images of both $q$ and $r$
lie in the boundary $\partial \oM_{g,n} \subset {\oM}_{g,n}$.

The cohomology and Chow groups 
of the moduli space of curves are
$$H^*({\oM}_{g,n},\mathbb{Q}) \ \ \ \text{and} \ \ \  A^*({\oM}_{g,n},\mathbb{Q})\, .$$
While there has been considerable progress in recent years,
many basic questions about the cohomology and
algebraic cycle theory
remain open.{\footnote{See
\cite{PandSLC} for a survey of results and open questions.}}

\subsection{Gromov-Witten classes}\label{gwcl}
Let $X$ be a nonsingular projective variety over $\mathbb{C}$, and let
$${\oM}_{g,n}(X,\beta)$$
be the moduli space of genus $g$, $n$-pointed 
stable maps to $X$ representing the class $\beta\in H_2(X,\mathbb{Z})$. The basic structures carried by
${\oM}_{g,n}(X,\beta)$
  are forgetful maps,
$$\pi: {\oM}_{g,n}(X,\beta) \rightarrow {\oM}_{g,n}\, ,$$
to the moduli space 
of curves via the domain
(in case $2g-2+n>0$) and
evaluation maps,
$$\text{ev}_i: {\oM}_{g,n}(X,\beta) \rightarrow X\, ,$$
for each marking $1\leq i \leq n$.

Given cohomology classes $v_1, \ldots, v_n \in H^*(X,\mathbb{Q})$,
the associated {\em Gromov-Witten class} is defined by
\begin{equation*}
\Omega_{g,n,\beta}^X(v_1, \ldots, v_n) = \pi_*\left(\prod_{i=1}^n \text{ev}_i^*(v_i) \, 
\cap\, [{\oM}_{g,n}(X,\beta)]^{\text{vir}}\right)\, \in\, H^*({\oM}_{g,n},\mathbb{Q})\, .
\end{equation*}
Central to the construction is 
the {\em virtual fundamental class} of the moduli space of stable maps,
$$[{\oM}_{g,n}(X,\beta)]^{\text{vir}}\in H_{2\cdot\text{virdim}}
({\oM}_{g,n}(X,\beta),\mathbb{Q})\, ,$$
of {\em virtual dimension} 
$$\text{virdim}= \int_{\beta}c_1(X) + (1-g)\cdot(\text{dim}_\C(X)-3) +n\,. $$ 
Gromov-Witten classes contain much more information than the
{\em Gromov-Witten invariants} defined by integration,
$$\big\langle v_1, \ldots, v_n\big\rangle^X_{g,n,\beta} = 
\int_{{\oM}_{g,n}}
\Omega_{g,n,\beta}^X(v_1, \ldots, v_n)\, 
.$$
We refer the reader to \cite{Beh,BehF,coxk,FulP}
for a detailed treatment of stable maps, virtual fundamental classes,
and Gromov-Witten invariants in algebraic geometry.

The Gromov-Witten classes satisfy formal properties
with respect to the natural  forgetful and boundary maps $p$, $q$, and $r$
discussed in Section \ref{xcxc}. 
The idea of a 
cohomological field theory was introduced by Kontsevich and Manin
\cite{KonMan} to fully capture these formal properties.

\subsection{Cohomological field theories}\label{cohax}
The starting point for defining a cohomological field theory is
a triple of data $(V,\eta,\b1)$ where
\begin{enumerate}
\item[$\bullet$]
$V$ is a finite dimensional $\mathbb{Q}$-vector space, 
\item[$\bullet$] $\eta$ is  
a non-degenerate symmetric 2-form on $V$,
\item[$\bullet$] $\b1 \in V$ is a distinguished element.
\end{enumerate}
Given a 
$\mathbb{Q}$-basis $\{e_i\}$ of $V$, the 
symmetric form $\eta$ can be written as a matrix
$$\eta_{jk}=\eta(e_j,e_k) \ .$$ The inverse matrix is denoted, as usual, by
$\eta^{jk}$.

A {\em cohomological field theory} consists of 
a system $\Omega = (\Omega_{g,n})_{2g-2+n > 0}$ of tensors
$$
\Omega_{g,n} \in H^*(\oM_{g,n},\mathbb{Q}) \otimes (V^*)^{\otimes n}.
$$
The tensor $\Omega_{g,n}$ associates a cohomology class in $H^*(\oM_{g,n},\mathbb{Q})$ 
to vectors 
$$v_1, \ldots, v_n\in V$$ assigned to the $n$ markings. We will use both 
$$\Omega_{g,n}(v_1\otimes \cdots \otimes v_n) \ \ \ \text{and} \ \ \
\Omega_{g,n}(v_1,\ldots, v_n) $$
to denote the associated cohomology class in $H^*(\oM_{g,n},\mathbb{Q})$.

In order to define a cohomological field theory, the system 
$\Omega = (\Omega_{g,n})_{2g-2+n > 0}$ must satisfy the 
CohFT axioms:

\begin{enumerate}
\item[(i)] Each tensor $\Omega_{g,n}$ is $\Sigma_n$-invariant for
 the natural action of the symmetric group $\Sigma_n$ on
$$H^*(\oM_{g,n},\mathbb{Q}) \otimes (V^*)^{\otimes n}$$
 obtained by simultaneously 
permuting the $n$ marked points of $\oM_{g,n}$ and the $n$ 
factors of $V^*$.
\item[(ii)] 
The tensor $q^*(\Omega_{g,n}) \in H^*(\oM_{g-1,n+2},\mathbb{Q}) \otimes 
(V^*)^{\otimes n}$,
obtained via pull-back by the boundary morphism 
$$
q : \oM_{g-1, n+2} \to \oM_{g,n}\ ,
$$
is required to equal the contraction
of $\Omega_{g-1,n+2}$ by the bi-vector 
$$\sum_{j,k} \eta^{jk} e_j \otimes e_k$$
inserted at the two identified points: 
$$q^*(\Omega_{g,n}(v_1,\ldots,v_n)) =
\sum_{j,k} \eta^{jk}\,\Omega_{g-1,n+2}(v_1,\ldots,v_n, e_j,e_k)$$
in $H^*(\oM_{g-1,n+2},\mathbb{Q})$ for all $v_i \in V$.

The tensor $r^*(\Omega_{g,n})$, obtained via pull-back 
by the boundary morphism
$$
r: \oM_{g_1, n_1+1} \times \oM_{g_2, n_2+1} \to \oM_{g,n}\, , 
$$
is similarly required
to equal the contraction of 
$\Omega_{g_1, n_1+1} \otimes \Omega_{g_2, n_2+1}$ by the
same bi-vector:
\begin{multline*}
r^*(\Omega_{g,n}(v_1,\ldots,v_n)) = \\
\sum_{j,k} \eta^{jk}\, \Omega_{g_1,n_1+1}(v_1,\ldots,v_{n_1}, e_j) \otimes
\Omega_{g_2,n_2+1}(v_{n_1+1}, \ldots, v_n, e_k)
\end{multline*}
in $H^*(\oM_{g_1,{n_1+1}},\mathbb{Q}) \otimes 
H^*(\oM_{g_2,{n_2+1}},\mathbb{Q})$ for all $v_i \in V$.

\item[(iii)] The tensor $p^*(\Omega_{g,n})$, obtained via pull-back by the forgetful map
$$p:\oM_{g,n+1} \rightarrow \oM_{g,n}\, ,$$
is required to satisfy 
$$
\Omega_{g,n+1}(v_1, \ldots,v_n , \b1) = p^*\Omega_{g,n} (v_1, \ldots,  v_n)\ 
$$
for all $v_i \in V$.
In addition, the equality 
$$\Omega_{0,3}(v_1,v_2,  \b1) = \eta(v_1,v_2)\, $$
is required for all $v_i \in V$.
\end{enumerate}

\begin{definition}\label{defcohft}
A system $\Omega= (\Omega_{g,n})_{2g-2+n>0}$ of tensors 
$$
\Omega_{g,n} \in H^*(\oM_{g,n},\mathbb{Q}) \otimes (V^*)^{\otimes n}
$$
satisfying (i) and (ii) is a {\em cohomological field theory} or a {\em CohFT}. If (iii) is also satisfied, $\Omega$ is 
a {\em CohFT with unit}.
\end{definition}

The simplest example of a cohomological field theory with unit is given by the
{\em trivial CohFT},
$$V=\mathbb{Q}\, ,\ \ \eta(1,1)=1\,, \ \ \b1=1 \,, \ \
\Omega_{g,n}(1,\ldots,1) = 1\in H^0(\oM_{g,n},\mathbb{Q})\, .$$
A more interesting example is given by the total Chern class 
$$c(\mathbb{E}) = 1 +\lambda_1+\ldots+ \lambda_g\, \in \, H^*(\oM_{g,n},\mathbb{Q})$$
of
the rank $g$ Hodge bundle $\mathbb{E}\rightarrow \oM_{g,n}$,
$$V=\mathbb{Q}\, ,\ \ \eta(1,1)=1 \,, \ \ \b1=1\,, \ \
\Omega_{g,n}(1,\ldots,1) = c(\mathbb{E})\in H^*(\oM_{g,n},\mathbb{Q})\, .$$

\begin{definition}\label{defcohft2}  For a CohFT
$\Omega= (\Omega_{g,n})_{2g-2+n>0}$, the {\em topological part}
$\omega$ of $\Omega$ is defined by 
$$\omega_{g,n} = [\Omega_{g,n}]^0 \, \in \, H^0(\oM_{g,n},\mathbb{Q}) \otimes
(V^*)^{\otimes n}
\, .$$
\end{definition}

The degree 0 part $[\,\, ]^0$ of $\Omega$ is simply obtained from the 
canonical summand projection
$$[\,\, ]^0:H^*(\oM_{g,n},\mathbb{Q}) \rightarrow H^0(\oM_{g,n},\mathbb{Q})\, .$$
If $\Omega$ is a CohFT with unit, then $\omega$ is also a
CohFT with unit.
The topological part of the CohFT obtained from the total
Chern class of the Hodge bundle is the trivial CohFT.

The motivating example of a CohFT with unit
is obtained from the Gromov-Witten
theory of a nonsingular projective variety $X$. Here,
$$V=H^*(X,\mathbb{Q})\, , \ \ \eta(v_1,v_2) = \int_X v_1\cup v_2\, , \ \ 
\b1 =  1\, .$$
Of course, the Poincar\'e pairing on $H^*(X,\mathbb{Q})$ 
 is symmetric only if $X$ has
no odd cohomology.{\footnote{To accommodate the case of arbitrary $X$,
the definition of a CohFT can be formulated with signs and
$\mathbb{Z}/2\mathbb{Z}$-gradings. We do not take the super vector space
path here.}} The
tensor $\Omega_{g,n}$ is defined using the Gromov-Witten classes $\Omega_{g,n,\beta}^X$
of Section \ref{gwcl} (together with a Novikov{\footnote{Formally, we
must extend scalars in the definition of a CohFT from $\mathbb{Q}$ to
the Novikov ring to capture the Gromov-Witten theory of $X$.}} parameter $q$),
$$\Omega_{g,n}(v_1,\ldots,v_n) = \sum_{\beta\in H_2(X,\mathbb{Q})}
\Omega_{g,n,\beta}^X\, q^\beta \,.$$
The CohFT axioms here coincide exactly with the axioms{\footnote{The divisor
axiom of Gromov-Witten theory (which concerns
divisor and curve classes on $X$) is not part of the CohFT axioms.}}
of Gromov-Witten theory related to the morphisms $p$, $q$, and $r$. 
For example, axiom (ii) of a CohFT
here is the {\em splitting
axiom} of Gromov-Witten theory, see \cite{KonMan}.

\subsection{Semisimplicity}

\label{sspty}
A CohFT  with unit $\Omega$ defines a {\em quantum product} $\bullet$ on $V$ 
by{\footnote{Since $\oM_{0,3}$ is a point, we canonically identify
$H^*(\oM_{0,3},\mathbb{Q}) \stackrel{\sim}{=} \mathbb{Q}$,
so $\Omega_{0,3}(v_1,v_2,v_3)\in \mathbb{Q}$.}}
$$\eta(v_1 \bullet v_2, v_3) = \Omega_{0,3}(v_1 \otimes v_2 \otimes v_3)\, .$$
The quantum product $\bullet$ is commutative by CohFT axiom (i). 
The associativity of $\bullet$ follows from CohFT axiom  (ii). The element
$\b1\in V$ is the identity for $\bullet$ by the second clause of
CohFT axiom (iii). Hence,
$$(V,\bullet, \b1)$$
is a commutative $\mathbb{Q}$-algebra.

\begin{lemma} \label{ggtt5} The topological part $\omega$ of $\Omega$ is uniquely and
effectively determined by the coefficients
$$\Omega_{0,3}(v_1,v_2,v_3) \in H^*(\oM_{0,3},\mathbb{Q})$$
of the quantum product $\bullet$.
\end{lemma}

\paragraph{Proof.}
Let the moduli point $[C,p_1,\ldots,p_n] \in\oM_{g,n}$ correspond to
 a maximally degenerate
curve (with every component isomorphic to $\mathbb{P}^1$ with exactly 3 special points).
Since 
$$\omega_{g,n}(v_1,\ldots,v_n) \in H^0(\oM_{g,n},\mathbb{Q})\, ,$$
is a multiple of the identity class, $\omega_{g,n}(v_1,\ldots,v_n)$
is determined by the pull-back to the point $[C,p_1,\ldots,p_n]$.
The equality
$$\omega_{g,n}(v_1,\ldots,v_n)\big|_{[C,p_1,\ldots,p_n]} =
\Omega_{g,n}(v_1,\ldots,v_n)\big|_{[C,p_1,\ldots,p_n]}\, $$
holds, and
the latter restriction is determined by 3-point values $\Omega_{0,3}(w_1,w_2,w_3)$
from repeated application of CohFT axiom (ii).
\qed

\vspace{10pt}

A finite dimensional $\mathbb{Q}$-algebra is {\em semisimple} if there
exists a basis $\{e_i\}$ of idempotents, 
$$e_ie_j = \delta_{ij} e_i\, ,$$
after an extension of scalars to $\com$.

\begin{definition}\label{defcohft3}
A CohFT with unit $\Omega= (\Omega_{g,n})_{2g-2+n>0}$ is {\em semisimple}
if $(V,\bullet, \b1)$ is a semisimple algebra.
\end{definition}

\subsection{Classification and calculation} \label{ww33}
The Givental-Teleman classification concerns semisimple CohFTs with unit.{\footnote{Semisimple CohFTs without unit are also covered, but we are
interested here in the unital case. Semisimplicity is an essential condition.}}
 The form of  the classification result
is as follows: {\em a semisimple CohFT with unit $\Omega$ is uniquely
determined by the following two structures:
\begin{enumerate}
\item[$\bullet$] the topological part $\omega$ of $\Omega$,
\item[$\bullet$] an $R$-matrix
$$R(z)= \mathsf{Id} + R_1 z + R_2 z^2 + R_3 z^3 + \ldots, \ \ \ R_k \in \text{\em End}(V)$$
satisfying the symplectic
property
$$R(z) \cdot R^\star(-z) = {\mathsf{Id}}\,, $$
where $\star$ denotes the adjoint with respect to the metric $\eta$.
\end{enumerate}}
\noindent The precise statement of the Givental-Teleman classification 
will be discussed in Section \ref{GTC}.

Via the Givental-Teleman classification, a semisimple CohFT  with unit 
$\Omega$ can be calculated in
three steps:
\begin{enumerate}
\item[(i)] determine the ring $(V,\bullet,\b1)$ as explicitly as possible,
\item[(ii)] find a closed formula for the topological part $\omega$ of
$\Omega$ via Lemma \ref{ggtt5},
\item[(iii)] calculate the $R$-matrix of the theory.
\end{enumerate}
In the language of Gromov-Witten theory, step (i) is the determination of
the {\em small quantum cohomology} ring $QH^*(X,\mathbb{Q})$ via 
the 3-pointed genus 0 Gromov-Witten invariants.
Step (ii) is then to calculate the Gromov-Witten invariants
where the domain has a {\em fixed} complex structure of higher
genus. 
New ideas are often required for the leap to higher genus moduli
in step (iii). Finding a closed formula for the $R$-matrix requires
a certain amount of luck.

Explaining how the above path to 
calculation plays out in three important CohFTs
 is my goal here. The three theories are:
\begin{enumerate}
\item[$\bullet$]
Witten's class on the moduli of $r$-spin curves, 
\item[$\bullet$] the Chern character of the
Verlinde bundle on the moduli of curves,
\item[$\bullet$]  the Gromov-Witten theory 
of the Hilbert scheme of points of $\com^2$. 
\end{enumerate}
While each theory has geometric interest and the
calculations have consequences in several directions, 
the focus of the paper will be on the CohFT determination.
The paths to calculation pursued here are applicable in many other cases.

\subsection{Past and future directions}
The roots of the classification of semisimple CohFTs can be found
in Givental's analysis \cite{Givental,Givental2,YPP} of the torus 
localization formula \cite{grpan} for the
higher genus Gromov-Witten theory of toric varieties. The 
three CohFTs treated here are not directly accessible via the older torus
localization methods. Givental's approach to the $R$-matrix via
oscillating integrals (used often in the study of toric geometries)
 is not covered in the paper.

Many interesting CohFTs are {\em not} semisimple. For example, 
the Gromov-Witten theory of the
famous Calabi-Yau quintic 3-fold, 
$$X_5 \subset \mathbb{P}^4\, ,$$
does {\em not} define a semisimple CohFT. However, in the past year,
an approach to the quintic via the semisimple {\em formal quintic} theory 
\cite{JGR,LhoP}
appears possible. These developments are not surveyed here.
 
\subsection{Acknowledgments}
Much of what I know about the Givental-Teleman classification
was learned through writing \cite{YPP} with Y.-P. Lee and
 \cite{PPZ} with A. Pixton and D. Zvonkine. 
For the study of the three CohFTs discussed in the paper, my collaborators
have been J. Bryan, 
F. Janda, A. Marian, A. Okounkov, D. Oprea, 
A. Pixton, H.-H. Tseng, and D. Zvonkine.
More specifically, 
\begin{enumerate}
\item[$\bullet$] Sections \ref{GTC}-\ref{wittenr} are based on 
the papers \cite{PPZ,PPZ2} and the Appendix of \cite{PPZ2},
\item[$\bullet$] Section \ref{verr} is based on the paper \cite{MOPPZ},
\item[$\bullet$] Section \ref{hhilb} is based on the papers \cite{bp,op} and especially
\cite{HHP}. 
\end{enumerate}
Discussions with A. Givental, T. Graber, H. Lho, and Y. Ruan
have played an important role in my view of the subject.
I was partially supported by
 SNF grant 200021-143274,   ERC grant
AdG-320368-MCSK, SwissMAP, and the Einstein Stiftung.

\section{Givental-Teleman classification}\label{GTC}

\subsection{Stable graphs} \label{stgr}
The boundary
strata of the moduli space of curves correspond
to {\em stable graphs} 
$$\Gamma=(\V, \H,\L, \ \mathrm{g}:\V \rarr \Z_{\geq 0},
\ v:\H\rarr \V, 
\ \mathrm{i} : \H\rarr \H)$$
satisfying the following properties:
\begin{enumerate}
\item[(i)] $\V$ is a vertex set with a genus function $\g:V\to \Z_{\geq 0}$,
\item[(ii)] $\H$ is a half-edge set equipped with a 
vertex assignment $v:H \to V$ and an involution $\mathrm{i}$,
\item[(iii)] $\E$, the edge set, is defined by the
2-cycles of $\mathrm{i}$ in $\H$ (self-edges at vertices
are permitted),
\item[(iv)] $\L$, the set of legs, is defined by the fixed points of $\mathrm{i}$ and endowed with a bijective correspondence with a set of markings,
\item[(v)] the pair $(\V,\E)$ defines a {\em connected} graph,
\item[(vi)] for each vertex $v$, the stability condition holds:
$$2\g(v)-2+ \n(v) >0,$$
where $\n(v)$ is the valence of $\Gamma$ at $v$ including 
both half-edges and legs.
\end{enumerate}
An automorphism of $\Gamma$ consists of automorphisms
of the sets $\V$ and $\H$ which leave invariant the
structures $\mathrm{g}$, $v$, and $\mathrm{i}$ (and hence respect $\E$ and $\L$).
Let $\text{Aut}(\Gamma)$ denote the automorphism group of $\Gamma$.

The genus of a stable graph $\Gamma$ is defined by:
$$g(\Gamma)= \sum_{v\in \V} \g(v) + h^1(\Gamma)\, .$$
Let $\mathsf{G}_{g,n}$ denote the set of all stable graphs
(up to isomorphism) of genus $g$ with $n$ legs. 
The strata{\footnote{We consider here the
standard stratification by topological type of the
pointed curve}} of the moduli space $\oM_{g,n}$ 
of Deligne-Mumford stable curves are in bijective correspondence to
 $\mathsf{G}_{g,n}$ by considering the dual graph of a generic pointed curve parameterized by the stratum.

To each stable graph $\Gamma$, we associate the moduli space
\begin{equation*}
\oM_\Gamma =\prod_{v\in \V} \oM_{\g(v),\n(v)}\, .
\end{equation*}
 Let $\pi_v$ denote the projection from $\oM_\Gamma$ to 
$\oM_{\g(v),\n(v)}$ associated to the vertex~$v$.  There is a
canonical
morphism 
\begin{equation}\label{dwwd}
\iota_{\Gamma}: \oM_{\Gamma} \rarr \oM_{g,n}
\end{equation}
 with image{\footnote{
The degree of $\iota_\Gamma$ is $|\text{Aut}(\Gamma)|$.}}
equal to the boundary stratum
associated to the graph $\Gamma$.  
%To construct $\xi_\Gamma$, 
%a family of stable pointed curves over $\oM_\Gamma$ is required.  Such a family
%is easily defined 
%by attaching the pull-backs of the universal families over each of the 
%$\oM_{\g(v),\n(v)}$  along the sections corresponding to half-edges.

\subsection {$R$-matrix action} \label{rmata}
\subsubsection{First action}
Let  $\Omega=(\Omega_{g, n})_{2g-2+n>0}$ be a CohFT{\footnote{$\Omega$ is not 
assumed here to be unital -- only CohFT axioms (i) and
(ii) are imposed.}}  
on the vector space $(V, \eta)$.
Let $R$ be a matrix series
$$R(z)= \sum_{k=0}^\infty R_k z^k \,   
\in\, {\mathsf{Id}} +z\cdot \text{End}(V)[[z]]$$ which satisfies the 
symplectic condition
$$ R(z) \cdot R^\star(-z) = \mathsf{Id}\, .$$
We define a
 new CohFT $R\Omega$
on the vector space $(V, \eta)$ by summing over stable graphs $\Gamma$
with summands given by  
 a product of vertex, edge, and leg contributions, 
\begin{equation}\label{romega}(R\Omega)_{g, n}=
\sum_{\Gamma\in\mathsf{G}_{g,n}} \frac{1}{|\text{Aut }(\Gamma)|} \iota_{\Gamma \star} \left(
\prod_{v\in \V}\mathsf {Cont}(v)\prod_{e\in \E}\mathsf {Cont }(e)
\prod_{l\in \L}\mathsf {Cont}(l)
 \right)\, ,\end{equation} where 
\begin {itemize}
\item [(i)] the vertex contribution is $$\mathsf {Cont}(v)=\Omega_{g(v), n(v)},$$ where $g(v)$ and $n(v)$ denote the genus and number of half-edges and legs of the vertex,
\item [(ii)] the leg contribution is the $\text{End}(V)$-valued cohomology class $$\mathsf{Cont}(l)=R(\psi_l)\, ,$$ where $\psi_l\in H^2(\oM_{g(v),n(v)},\mathbb{Q})$ 
is the cotangent class at the marking corresponding to the leg, 
\item [(iii)] the edge contribution is $${\mathsf {Cont}}(e)=\frac{\eta^{-1}-R(\psi'_e)\eta^{-1} R(\psi''_e)^{\top}}{\psi'_e+\psi''_e}\, , $$ 
where $\psi'_e$ and $\psi''_e$ are the cotangent classes at the node which represents the edge $e$. The symplectic condition guarantees that the edge contribution is well-defined. 
\end{itemize}

We clarify the meaning of the edge contribution (iii), 
$$
\mathsf{Cont}(e)\in V^{\otimes 2}\otimes H^{\star}(\oM_{g', n'})\otimes H^{\star}(\oM_{g'', n''}),
$$ 
where $(g', n')$ and $(g'', n'')$ are the labels of the vertices adjacent to
the edge $e$ by writing the formula explicitly in coordinates.

Let $\{e_\mu\}$ be a $\mathbb{Q}$-basis of $V$. The components of the $R$-matrix in the basis are $R_\mu^\nu(z)$,
$$
R(z)(e_\mu)=\sum_{\nu} R_\mu^\nu (z) \cdot e_\nu\, .
$$ 
The components of $\mathsf{Cont}(e)$ are 
$$
\mathsf{Cont}(e)^{\mu \nu} = 
\frac{\eta^{\mu \nu}- \sum_{\rho, \sigma} R_\rho^\mu(\psi_e')\cdot \eta^{\rho \sigma}\cdot R_\sigma^\nu(\psi_e'')}{\psi_e'+\psi_e''}\in H^{\star}(\oM_{g', n'})\otimes H^{\star}(\oM_{g'', n''}).
$$ 
The fraction 
$$
\frac{\eta^{\mu \nu}-\sum_{\rho, \sigma}R_\rho^\mu(z)\cdot \eta^{\rho \sigma}\cdot R_\sigma^\nu(w)}{z+w}
$$ 
is a power series in $z$ and $w$ since the numerator vanishes when $z=-w$
as a consequence of the symplectic condition which, in coordinates,  takes
 the form 
$$
\sum_{\rho, \sigma} R_\rho^\mu(z)\cdot \eta^{\rho \sigma}\cdot R_\sigma^\nu(-z) = \eta^{\mu \nu}.
$$ 
The substitution $z=\psi_e'$ and $w=\psi_e''$ is therefore unambiguously defined.

\begin{definition}\label{tttt}
Let $R\Omega$ be the CohFT obtained from $\Omega$
by the $R$-action \eqref{romega}.
\end{definition}

 The above 
$R$-action was first defined{\footnote{To simplify our formulas, we have changed Givental's and Teleman's conventions by replacing $R$ with $R^{-1}$. 
Equation \eqref{romega} above then determines a right group action on CohFTs rather than a left group action as in Givental's and Teleman's papers.}}\,on Gromov-Witten potentials by Givental~\cite{Givental}. 
An abbreviated treatment of the lift to CohFTs appears in
papers by Teleman \cite{Teleman} and Shadrin \cite{Shadrin}.
A careful proof that $R\Omega$ satisfies  CohFT axioms (i) and (ii)
 can be found in
\cite[Section 2]{PPZ}. 

If $\Omega$ is a CohFT with unit on $(V,\eta,\b1)$, 
then $R\Omega$
may not respect the unit $\b1$. To handle the unit, a second action is
required.

\subsubsection{Second action}

A second action on the CohFT{\footnote{To define
the translation action, $\Omega$
is required only to be CohFT and not necessarily a CohFT with unit.}}  
$\Omega$ on $(V,\eta)$
is given by translations. Let 
$T\in V[[z]]$ 
be a series with {\em no terms of degree $0$ or $1$},
$$T(z)=T_2z^2+T_3z^3+\ldots\, ,\ \ \ T_k\in V\, .$$

\begin{definition}\label{tttt5}
Let $T\Omega$ be the CohFT obtained from $\Omega$ by the formula
\begin{equation*}%\label{translation} 
(T\Omega)_{g, n} (v_1, \ldots, v_n)=\sum_{m=0}^{\infty} \frac{1}{m!} p_{m\star} 
\Big(\Omega_{g, n+m}(v_1, \ldots, v_n, T(\psi_{n+1}),\ldots, T(\psi_{n+m}))\Big)
\,, 
\end{equation*} 
where $p_{m}:\oM_{g, n+m}\to \oM_{g, n}$ is the morphism forgetting
the last $m$ markings.
\end{definition}

The right side of the formula
in Definition \ref{tttt5} is a formal expansion
by distributing the powers of the $\psi$ classes as follows: 
$$\Omega_{g, n+ m}(\cdots, T(\psi_{\bullet}),\cdots )=\sum_{k=2}^{\infty} \psi_{\bullet}^k \cdot \Omega_{g, n+m} (\cdots, T_k , \cdots)\, .$$
 The summation is finite because $T$ has no terms of degree 0 or 1.

 \subsection {Reconstruction} 
We can now state the Givental-Teleman classification result \cite{Teleman}. 
Let $\Omega$ be a semisimple CohFT with unit on $(V, \eta,\b1)$, and
let $\omega$ be the topological part of $\Omega$. 
For a symplectic matrix $R$, define 
$$R.\omega = R(T(\omega)) \ \ \text{with} \ \ 
T(z)=z(\left({\mathsf{Id}}-R(z)\right)\cdot \b1)\in V[[z]]\, .
$$
By \cite[Proposition 2.12]{PPZ}, $R.\omega$ is a CohFT with unit on $(V,\eta, \b1)$.
The Givental-Teleman classification asserts
the {\em existence} of a unique  $R$-matrix which exactly
recovers $\Omega$. 

\begin{theorem} \label{GTc}
There exists a unique symplectic matrix 
$$
R\in {\mathsf{Id}}+z\cdot \End (V)[[z]]
$$
which reconstructs $\Omega$ from $\omega$,
$$
\Omega=R.\omega\,,
$$
as a CohFT with unit.
%where the series $T$ is defined 
%by
%$$
%T(z)=z(\left({\mathsf{Id}}-R(z)\right)\cdot \b1)\in V[[z]]\, .
%$$ 
\end{theorem}

%A proof of the uniqueness of $R$
%can be found in \cite{}. 
The first example concerns
the total Chern class CohFT of Section \ref{cohax},
$$V=\mathbb{Q}\, ,\ \ \eta(1,1)=1\,, \ \ \b1=1 \,, \ \
\Omega_{g,n}(1,\ldots,1) = c(\mathbb{E})\in H^*(\oM_{g,n},\mathbb{Q})\, .$$
The topological part is the trivial CohFT, and the
$R$-matrix is
$$R(z)= \exp\left(-\sum_{k=1}^\infty \frac{B_{2k}}{(2k)(2k-1)} z^{2k-1}\right)\, .$$
That the above $R$-matrix reconstructs the total Chern class CohFT
is a consequence  of Mumford's calculation \cite{Mum}
of the Chern character  of the Hodge bundle by Grothendieck-Riemann-Roch.

\subsection{Chow field theories}
Let $(V,\eta, \b1)$ be a $\mathbb{Q}$-vector space with
a non-degenerate symmetric 2-form and a distinguished element. 
Let $\Omega=(\Omega_{g,n})_{2g-2+n>0}$ be 
a system of tensors
$$
\Omega_{g,n} \in A^*(\oM_{g,n},\mathbb{Q}) \otimes (V^*)^{\otimes n}
$$
where $A^*$ is the Chow group of algebraic cycles modulo
rational equivalence.
In order to define a Chow field theory, the system 
$\Omega$ must satisfy the 
CohFT axioms of Section \ref{cohax}  with cohomology $H^*$ replaced everywhere 
by Chow $A^*$.

\begin{definition}\label{defcohft9}
A system $\Omega= (\Omega_{g,n})_{2g-2+n>0}$ of elements 
$$
\Omega_{g,n} \in A^*(\oM_{g,n},\mathbb{Q}) \otimes (V^*)^{\otimes n}
$$
satisfying (i) and (ii) is a {\em Chow field theory} or a {\em ChowFT}. If (iii) is also satisfied, $\Omega$ is 
a {\em ChowFT with unit}.
\end{definition}

For ChowFTs,
the quantum product $(V,\bullet, \b1)$ and semisimplicity are defined
just as  for CohFTs. The $R$- and $T$-actions of Sections \ref{rmata}
also lift immediately to ChowFTs. However,
the classification of semisimple ChowFTs is
an open question.

\begin{question}\label{ffxx2} Does the Givental-Teleman classification of Theorem \ref{GTc} 
hold for a semisimple Chow field theory $\Omega$ with unit?
\end{question}

\section{Witten's $r$-spin class}
\label{wittenr}

\subsection{$r$-spin CohFT} \label{wsc}
Let $r\geq 2$ be an integer. Let $(V_r,\eta,\b1)$ 
be the following triple: 
\begin{enumerate}
\item[$\bullet$] $V_r$ is an $(r-1)$-dimensional $\mathbb{Q}$-vector space with basis
 $e_0, \dots, e_{r-2}$, 
\item[$\bullet$] $\eta$ is the non-degenerate symmetric $2$-form
$$
\eta_{ab} = \langle e_a, e_b \rangle =\delta_{a+b,r-2} \, ,
$$
\item[$\bullet$] $\b1 = e_0$. 
\end{enumerate}
Witten's $r$-spin theory provides a family of classes
$$
\mathcal{W}^r_{g,n}(a_1, \dots, a_n) \in H^*(\oM_{g,n},\mathbb{Q})\, 
$$ 
for $a_1, \dots, a_n \in \{0, \dots, r-2 \}$
which define a CohFT $\W^r=(\W^r_{g,n})_{2g-2+n>0}$ by 
$$
\W^r_{g,n}: V^{\otimes n}_r \rightarrow H^*(\oM_{g,n},\mathbb{Q})\, , 
\ \ \ \
\W^r_{g,n}( e_{a_1} \otimes \cdots \otimes e_{a_n}) =
\mathcal{W}^r_{g,n}(a_1, \dots, a_n)\, .
$$
%To emphasize $r$, we will often refer to $V$ as $V_r$.
The class $\mathcal{W}^r_{g,n}(a_1, \dots, a_n)$ has (complex) 
degree{\footnote{So  $\mathcal{W}^r_{g,n}(a_1, \dots, a_n) \in H^{2\cdot \D^r_{g,n}(a_1,\ldots,a_n)}(\oM_{g,n},\mathbb{Q})$.}}
\begin{eqnarray}
\label{gred}
\text{deg}_{\C}\ \mathcal{W}^r_{g,n}(a_1, \dots, a_n) & = & 
\D^r_{g,n}(a_1, \dots, a_n) \\ \nonumber
& = & \frac{(r-2)(g-1) + \sum_{i=1}^n a_i}{r}\ .
\end{eqnarray}
If $\D^r_{g,n}(a_1, \dots, a_n)$ is not an integer, the corresponding
Witten's class vanishes.

The construction of $\mathcal{W}_{0,n}^r(a_1,\ldots,a_n)$
in genus 0 was carried out by Witten \cite{Witten}
using $r$-spin structures. 
Let $\oM_{0,n}^r(a_1,\ldots,a_n)$ be the Deligne-Mumford moduli space
parameterizing  $r^{\rm th}$ roots, 
$$\mathcal{L}^{\otimes r} \stackrel{\sim}{=} \omega_C\Big(-\sum_{i=1}^n a_i p_i\Big) \ \ \ \text{where} \ \ \ [C,p_1,\ldots,p_n]\in
\oM_{0,n} \, .$$
The class $\frac{1}{r}\mathcal{W}^r_{0,n}(a_1,\ldots,a_n)$ is defined to be the
push-forward to $\oM_{0,n}$ of the top Chern class of the 
bundle on $\oM^{r}_{0,n}(a_1,\ldots,a_n)$
 with fiber
$H^1(C,\mathcal{L})^*$.

The existence of Witten's class in higher genus is both remarkable and
highly non-trivial. 
Polishchuk and Vaintrob
\cite{Pol,PolVai} constructed
$$\mathcal{W}^r_{g,n}(a_1, \dots, a_n)\in A^*(\oM_{g,n},\mathbb{Q})$$
as an algebraic cycle class and proved 
and the CohFT axioms (i-iii) for a Chow field theory hold.
 The algebraic approach was 
later simplified in \cite{ChLi,Chiodo}. 
Analytic constructions appear in \cite{fjr,Mochizuki}.

%and satisfies

%Whether Witten's $r$-spin theory as an algebraic cycle
%takes values in
%$$R^*(\cM_{g,n})\subset A^*(\cM_{g,n})$$
% is an open question.

\subsection{Genus 0} \label{gen00}
\subsubsection{3 and 4 markings}
Witten \cite{Witten} determined  the following initial conditions
in genus 0 with $n=3,4$:
\begin{equation}\label{fred}
\int_{\oM_{0,3}}\mathcal{W}^r_{0,3}(a_1,a_2,a_3) = 
\left|
\begin{array}{cl}
1 & \mbox{ if } a_1+a_2+a_3 = r-2\, ,\\
0 & \mbox{ otherwise,}
\end{array}
\right.
\end{equation}
$$
\int_{\oM_{0,4}}\mathcal{W}^r_{0,4}(1,1,r-2,r-2) = \frac1{r}\, .
$$
Uniqueness of the $r$-spin CohFT in
genus 0 follows easily from the initial conditions \eqref{fred}
and the axioms of a CohFT with unit.

The genus $0$  sector  of the CohFT $\W^r$ defines a quantum product{\footnote{See Section \ref{sspty}.}} $\bullet$
on $V_r$.
The resulting algebra $(V_r,\bullet,\b1)$, even after
extension to $\C$, is {\em not} semisimple.
Therefore, the Givental-Teleman classification can not be directly
applied.

\subsubsection{Witten's $r$-spin class and representations of $\slc$.}
Consider the Lie algebra $\sl_2=\slc$.
Denote by $\rho_k$ the $k^{th}$ symmetric power of the standard 2-dimensional representation of $\sl_2$,
$$\rho_k =\text{Sym}^k(\rho_1)\, , \ \ \   \dim \rho_k = k+1\, .$$
The following complete solution of the genus 0 part of the CohFT $\W^r$ 
(after integration) was
found by Pixton, see \cite{PPZ2} for a proof.

\begin{theorem} \label{Thm:sl2}
Let  ${\bf{a}}=(a_1, \dots, a_{n\geq 3})$ with 
$a_i \in \{ 0, \dots, r-2 \}$  satisfy the degree constraint
$\D^r_{0,n}({\bf{a}})= n-3$.
Then, 
\begin{equation*}
\int_{\oM_{0,n}}
\mathcal{W}^r_{0,n}({\bf{a}}) = 
\frac{(n-3)!}{r^{n-3}} 
\dim \Bigl[ \rho_{r-2-a_1} \otimes \cdots \otimes \rho_{r-2-a_n} \Bigr]^{\sl_2} ,
\end{equation*}
where the superscript $\sl_2$ denotes the $\sl_2$-invariant subspace.
\end{theorem}

%\vspace{5pt}
%The degree constraint $\D^r_{0,n}({\bf{a}})= n-3$ in the statement of
% Theorem \ref{Thm:sl2} can be written
% equivalently (using \eqref{gred}) as
%$$\sum_{i=1}^n a_i = (n-2)r-2\, .$$
%Since $a_i\leq r-2$, the bound $n\leq r+1$ is a simple consequence. 

\subsubsection{Shifted Witten class}
\label{frrz}

\begin{definition}\label{frrf}
For $\gamma \in V_r$, the {\em shifted} $r$-spin CohFT 
$\W^{r,\gamma}$
is defined by 
\begin{equation*}
\W_{g,n}^{r,\gamma}(v_1 \otimes \cdots \otimes v_n) = 
\sum_{m \geq 0}\frac{1}{m!}\,
p_{m\star}\, \W^r_{g,n+m}(v_1 \otimes \cdots \otimes v_n \otimes \gamma^{\otimes m}),
\end{equation*}
where $p_m \colon \oM_{g,n+m} \to \oM_{g,n}$ is the map forgetting
the last $m$ markings.
\end{definition}
Using degree formula $\D^r_{g,n}$, the summation in the definition
of the shift
is easily seen to be finite. 
The shifted Witten class $\W^{r,\gamma}$ determines a CohFT with unit, see
\cite[Section 1.1]{PPZ}.

\begin{definition}\label{hathat}
Define the CohFT $\widehat{\W}^r$ with unit on $(V_r,\eta,\b1)$ by the
shift 
$$\widehat{\W}^r=\W^{r,(0,\ldots,0,r)}$$
along the special vector $re_{r-2}\in V_r$.
Let $(V_r, \widehat{\bullet},\b1)$ be the $\mathbb{Q}$-algebra
determined by the quantum product defined by $\widehat{\W}^r$.
\end{definition}

The {\em Verlinde algebra} of level~$r$ for $\sl_2$ is spanned by the weights of $\sl_2$ from $0$ to $r-2$. The coefficient of $c$ in the product $a \bullet b$ is equal to the dimension of the $\sl_2$-invariant subspace of the 
representation $\rho_a \otimes \rho_b \otimes \rho_c$ provided the inequality 
$$a+b+c \leq 2r-4$$ is satisfied. 
Using Theorem \ref{Thm:sl2}  for the integral $r$-spin theory in genus 0, the
following basic result is proven in \cite{PPZ2}.

\begin{theorem}  \label{Thm:Verlinde}
The algebra 
$(V_r,\widehat{\bullet},\b1)$ is isomorphic to the Verlinde algebra of level~$r$ for $\sl_2$. 
\end{theorem}

Since the Verlinde algebra is well-known to be semisimple{\footnote{An explicit normalized
idempotent basis
is given in Proposition \ref{Prop:idempotents} below.}}, the Givental-Teleman
classification of Theorem \ref{GTc} can 
be applied to the CohFT $\widehat{\W}^r$.
Using the degree formula \eqref{gred} and Definition \ref{hathat}, we
see the (complex) degree $\D^r_{g,n}(a_1,\ldots,a_n)$ part of $\widehat{\W}^r$ equals
$\W^r$,
$$ \Big[ \widehat{\W}^r_{g,n}(a_1,\ldots,a_n)\Big]^{\D^r_{g,n}(a_1,\ldots,a_n)} = \W^r_{g,n}(a_1,\ldots,a_n)\, .$$
Hence, a complete computation of $\widehat{\W}^r$ also provides a
computation of $\W^r$.

\subsection{The topological field theory} \label{Ssec:TopFT12}
After the studying genus 0 theory, we turn our attention to the 
topological part $\widehat{\omega}^r$ of $\widehat{\W}^r$.
The following two results of \cite{PPZ2} provide a complete
calculation.% of $\widehat{\omega}^r$.

\begin{proposition} \label{Prop:idempotents}
The basis of normalized idempotents of $(V_r,\widehat{\bullet},\b1)$ is given by
$$
v_k = \sqrt{\frac2r} \; \sum_{a=0}^{r-2} \sin \left(
\frac{(a+1)k\pi}r \right) \; e_a, \qquad k\in\{ 1, \dots, r-1\}\, .
$$
More precisely, we have
$$
\eta(v_k, v_l) = (-1)^{k-1} \delta_{k,l}\, ,
\qquad \ \
v_k \,\widehat{\bullet}\, v_l = %\phi^{(r-2)/4} 
\frac{\sqrt{r/2}}{\sin(\frac{k \pi}{r})} v_k \, \delta_{k,l}\, .
$$
\end{proposition}

Once the normalized idempotents are found, the computation of
 $\widehat{\omega}^r$ is
straightforward by Lemma \ref{ggtt5} and elementary trigonometric
identities.

\begin{proposition} \label{Prop:QuantumProd1}
For $a_1, \dots, a_n \in \{0, \dots, r-2 \}$,
we have
\begin{equation*}
\widehat{\omega}^{r}_{g,n}(e_{a_1}, \ldots,e_{a_n}) =
\left(\frac{r}2\right)^{g-1}
 \; 
\sum_{k=1}^{r-1}
\frac
{(-1)^{(k-1)(g-1)}\prod\limits_{i=1}^n \sin \left(
\frac{(a_i+1)k\pi}r \right)}
{\left(\sin(\frac{k\pi}{r})\right)^{2g-2+n}}\, .
\end{equation*}
\end{proposition}

\noindent In Proposition \ref{Prop:QuantumProd1}, the CohFT $\widehat{\omega}^r$
is viewed as taking values in $\mathbb{Q}$ via the
canonical identification $$H^0(\oM_{g,n},\mathbb{Q}) \stackrel{\sim}{=} \mathbb{Q}\, . $$
\subsection{The $R$-matrix} \label{Ssec:Rmatrix1}
The last (and often hardest) 
step in the computation of a semisimple CohFT 
via the Givental-Teleman classification is to find the unique $R$-matrix.
Remarkably, there exists a closed formula in hypergeometric series for
the $R$-matrix of the CohFT $\widehat{\W}^r$.
The precise shift in Definition \ref{hathat} of the CohFT $\widehat{\W}^r$
is crucial: the shift $re_{r-2}$ is (up to scale) the
{\em only} shift of $\W^r$ for which 
closed formulas for the $R$-matrix are known. 

The method of finding the unique $R$-matrix for 
$\W^r$ uses the {\em Euler field} $e_{r-2}$ at the shift $re_{r-2}$.
The operator of quantum multiplication $\widehat{\bullet}$ 
by the Euler field 
in the basis $e_0,\ldots, e_{r-2}$ is 
$$
\xi = 
{\small{\begin{pmatrix}
\nice 0 & \nice \cdots & \nice \cdots & \nice 0  & \nice 2\\
 \qquad  & \qquad & \qquad & \qquad & \qquad \\
\nice 0 & & \nice 0 & \nice 2 & 0 \\
\\
\nice \vdots & \nice \iddots & \nice \iddots  & \nice \iddots & \nice \vdots \\
\\
\nice 0 & \nice 2 & \nice 0  &  & \nice 0 \\
\\
\nice 2 & \nice 0 & \nice \cdots  & \nice \cdots & \nice 0 
\end{pmatrix}\,}} .
$$
In the same frame, the {\em shifted degree operator} is
$$
\mu = \frac1{2r} 
{\small{\begin{pmatrix}
\nice -(r-2) & \nice 0 & \nice \cdots & \nice \cdots & \nice 0  \\
 \qquad & \qquad  & \qquad  & \qquad & \qquad \\
\nice 0 & -(r-4) & \nice 0 &  & 0 \\
\\
\nice \vdots & \nice \ddots & \nice \ddots  & \nice \ddots & \nice \vdots \\
\\
\nice 0 &  & \nice 0  & r-4 & \nice 0 \\
\\
\nice 0  & \nice \cdots  & \nice \cdots & \nice 0 & \nice r-2 
\end{pmatrix}\,}} .
$$
Since $\widehat{\W}^r$ has an Euler field with an associated degree operator, the
unique $R$-matrix for the classification is given by the solution of 
\begin{equation}\label{ff4455}
[R_{m+1},\xi] = (m-\mu) R_m
\end{equation}
with the initial condition $R_0=\mathsf{Id}$, see \cite{Teleman}.

\begin{definition} 
For each integer $a \in \{0, \dots, r-2 \}$, define the hypergeometric series 
$$
\B_{r,a}(z) = \sum_{m=0}^\infty
\left[
\prod_{i=1}^m 
\frac{\big((2i-1)r-2(a+1)\big)\big((2i-1)r+2(a+1)\big)}i 
\right] \!\! 
\left( - \frac{z}{16r^2} \right)^{\! m} \!\!\!.
$$
Let $\Be_{r,a}$ and $\Bo_{r,a}$ the even and odd summands{\footnote{The even
summand consists of all the even powers of $z$ (and likewise for the
odd summand).}} of the 
series~$\B_{r,a}$.
\end{definition}

%The unique solution $R(z)=\sum_{m=0}^\infty R_m z^m \in \End(V_r)[[z]]$ of the equations
%$$
%[R_{m+1},\xi] = (m+\mu) R_m
%$$
%with the initial condition $R_0=1$ has coefficients
%$$
%R^a_a = \Be_{r,r-2-a}(\phi^{-r/2}z)\,, \ \ \ \  a \in \{0, \dots, r-2 \}
%$$
%on the main diagonal,
%$$
%R^{r-2-a}_a =  -\Bo_{r,a}(\phi^{-r/2}z)\,, \ \ \ \  
%a \in \{0, \dots, r-2 \}
%$$
%on the antidiagonal (if $r$ is even, the coefficient at the intersection of both diagonals is 1), and $0$ e%verywhere else. 

The unique solution to \eqref{ff4455} is computed in \cite{PPZ2}. The $R$-matrix
of $\widehat{\W}^r$ has a surprisingly simple form.

\begin{theorem}
The unique $R$-matrix classifying $\widehat{\W}^r$ has coefficients
$$
R^a_a = \Be_{r,a}(z)\,, \ \ \ \ 
 a \in \{0, \dots, r-2 \}
$$
on the main diagonal, and
$$
R^{r-2-a}_a = \Bo_{r,a}(z)\,, \ \ \ \ 
 a \in \{0, \dots, r-2 \}
$$
on the antidiagonal (if $r$ is even, the coefficient at the intersection of both diagonals is 1), and $0$ everywhere else. 
\end{theorem}

In case $r=2$, the matrix is trivial $R(z)=\mathsf{Id}$. For 
$r=3$ and $4$ respectively, the $R$-matrices{\footnote{$R$ here is
$R^{-1}$ in \cite{PPZ,PPZ2} because of a change of conventions.}} are
$$R(z)=
\begin{pmatrix}
\Be_{3,0}(z) & \Bo_{3,1}(z) \\
\Bo_{3,0}(z) & \Be_{3,1}(z)
\end{pmatrix}\, , $$
$$
R(z)=\begin{pmatrix}
\ \ \ \Be_{4,0}(z) & \ 0\  & \Bo_{4,2}(z)   \\
0 & \ 1\  & 0 \\
\Bo_{4,0}(z) &\ 0\  & \ \,  \Be_{4,2}(z)
\end{pmatrix}\, . $$

\subsection{Calculation of $\W^r$}
The analysis of Sections \ref{gen00}-\ref{Ssec:Rmatrix1} together complete the
calculation of $\widehat{\W}^r$,
$$\widehat{\W}^r = R.\widehat{\omega}^r\,,$$
in exactly the steps (i)-(iii) proposed in Section \ref{ww33} of the Introduction.
Then, as we have seen,
$$\W^r_{g,n}(a_1,\ldots,a_n) = \Big[ \widehat{\W}^r_{g,n}(a_1,\ldots,a_n)\Big]^{\D^r_{g,n}(a_1,\ldots,a_n)}\, .$$
The calculation has an immediate consequence \cite{PPZ2}.

\begin{corollary} \label{ddrr14} Witten's $r$-spin class on $\oM_{g,n}$ lies in the tautological ring (in cohomology),
$$\mathcal{W}^r_{g,n}(a_1,\ldots,a_n) \in RH^*(\oM_{g,n},\mathbb{Q})\, .$$
\end{corollary}

We refer the reader to \cite{PandSLC} for a discussion of tautological classes
on the moduli space of curves. In fact, the first proof of Pixton's relations
in $RH^*(\oM_{g,n},\mathbb{Q})$  was
obtained via 
 the calculation
of $\widehat{\W}^3$ in \cite{PPZ}.

\subsection{Questions}
Whether Corollary \ref{ddrr14} also holds in Chow is an interesting question: is
\begin{equation}\label{qqqq}
\mathcal{W}^r_{g,n}(a_1,\ldots,a_n) \in R^*(\oM_{g,n},\mathbb{Q})\, ?
\end{equation}
A positive answer to Question \ref{ffxx2} about the classification of Chow field
theories would imply a positive answer here. The following question may be
viewed as a  refinement of  \eqref{qqqq}.

\begin{question}  Find a formula in algebraic cycles
for Witten's r-spin class on  $\oM_{g,n}^r(a_1,\ldots,a_n)$
{\em before} push-forward to $\oM_{g,n}$.
\end{question}

Another open direction concerns 
the moduli spaces of holomorphic differentials \cite{B5,FarP}.
Let $(a_1,\ldots,a_n)$ be a partition of $2g-2$
with non-negative parts. Let
\vspace{-3pt}
$$\overline{\HH}_{g}(a_1, \dots, a_n)\subset \oM_{g,n}\, $$

\vspace{-2pt}
\noindent be the closure of the locus of moduli points
\vspace{-10pt}
$$[C,p_1,\ldots,p_n]\in \cM_{g,n} \ \ \ \text{where} \ \ \
\omega_C \stackrel{\sim}{=} \mathcal{O}_C\Big(\sum_{i=1}^n a_i p_i\Big)\, .$$

\vspace{-4pt}
\noindent For 
$r-2\geq \text{Max}\{a_1,\ldots,a_n\}$,
Witten's $r$-spin class 
$\mathcal{W}^r_{g,n}(a_1, \dots, a_n)$
is well-defined and of degree {\em independent} of $r$,
\vspace{-2pt}
$${\mathsf{D}}^r_{g,n}(a_1,\ldots,a_n) = \frac{(r-2)(g-1) + \sum_{i=1}^n a_i}{r} =g-1\, .$$
By \cite[Theorem 7]{PPZ2},
 after scaling by $r^{g-1}$, 
\vspace{-2pt}
$$
\mathcal{W}_{g,n}(a_1, \dots, a_n)[r]=r^{g-1}\cdot \mathcal{W}^r_{g,n}(a_1, \dots, a_n)\in RH^{g-1}(\oM_{g,n},\mathbb{Q})
$$
is a {\em polynomial} in $r$ for all sufficiently large $r$.

\begin{question}
Prove the following conjecture of \cite[Appendix]{PPZ2}:
$$
(-1)^g \mathcal{W}_{g,n}(a_1, \dots, a_n)[0] = [\overline{\HH}_{g}(a_1, \dots, a_n)]\in H^{2(g-1)}(\oM_{g,n},\mathbb{Q})\, .
$$
\end{question}

\section{Chern character of the Verlinde bundle} \label{verr}

\subsection {Verlinde CohFT} Let $G$ be a complex, simple, 
simply connected Lie group with Lie algebra $\mathfrak g$.
Fix an integer {\em level} $\ell>0$. Let $(V_\ell,\eta, \b1)$ be the
following triple:
\begin{enumerate}
\item[$\bullet$] $V_\ell$ is the $\mathbb{Q}$-vector space
with basis indexed by the irreducible representations of $\mathfrak g$ at level $\ell$,
\item[$\bullet$] $\eta$ is the non-degenerate symmetric 2-form
 $$\eta(\mu, \nu)=\delta_{\mu, \nu^{\star}}$$ where $\nu^{\star}$ denotes the dual representation,
\item[$\bullet$] $\b1$ is the basis element corresponding
 to the trivial representation. 
\end{enumerate}

Let   
 $\mu_1, \dots, \mu_n$
be $n$ irreducible representations of  $\mathfrak g$ at level $\ell$.
A vector bundle $${\mathbb V}_{g}(\mu_1, \dots, \mu_n)\to \oM_{g,n}$$ 
is constructed in \cite {TUY}. Over nonsingular curves, the fibers of 
${\mathbb V}_{g}(\mu_1, \dots, \mu_n)$ are the spaces of non-abelian theta functions -- spaces of global sections of the determinant line bundles over the moduli of parabolic $G$-bundles. 
To extend ${\mathbb V}_{g}(\mu_1, \dots, \mu_n)$ over 
the boundary $$\partial\cM_{g,n}\subset \oM_{g,n}\,,$$ the
theory of {\it conformal blocks} is required \cite {TUY}.
The vector bundle
 $\mathbb V_g(\mu_1, \ldots, \mu_n)$ has various names in
the literature:
the Verlinde bundle, the bundle of conformal blocks, and the bundle of vacua.
A study in genus 0 and 1 can be found in \cite{F}.

A CohFT $\Omega^\ell$ is defined via the Chern
 character{\footnote{For a vector bundle $\mathbb V$ with Chern roots $r_1, \ldots, r_k$, $$\ch_t(\mathbb V)=\sum_{j=1}^{k} e^{tr_j}\, .$$ 
The parameter $t$ may be treated either as a formal variable, in which case 
the CohFT is defined over the ring $\mathbb{Q}[[t]]$ instead of 
$\mathbb{Q}$, or as a rational number $t \in \mathbb{Q}$.}}
 of the
Verlinde bundle: 
$$\Omega^\ell_{g, n}(\mu_1, \ldots, \mu_{n})=\ch_t(\mathbb V_{g}(\mu_1, \ldots, \mu_n))\in H^{\star}(\oM_{g, n},\mathbb{Q})\, .$$ 
CohFT axiom (i) for $\Omega^\ell$  is trivial. 
Axiom (ii) follows from the fusion rules \cite {TUY}.
Axiom (iii) for the unit $\b1$
 is the {\em propagation of vacua} \cite[Proposition 2.4(i)]{F}.

\subsection{Genus 0 and the topological part}
%The CohFT $\Omega^\ell$ is semisimple. 
Since the variable $t$ carries the degree grading, the topological
part $\omega^\ell$ of $\Omega^\ell$ is
obtained by setting $t=0$, 
$$
\omega^\ell_{g,n}=\Omega^\ell_{g, n}\Big|_{t=0}\, .
$$ 
The result is the just the rank of the Verlinde bundle,
$$
\omega^\ell_{g, n}(\mu_1, \ldots, \mu_n)=\text{rk } {\mathbb V}_g(\mu_1, \ldots, \mu_n)=d_g(\mu_1, \ldots, \mu_n)\, .
$$ 
With the quantum product obtained{\footnote{Since
the quantum product depends only upon the tensors of genus 0
with 3 markings, the quantum products of $\Omega^\ell$ and
$\omega^\ell$ are equal.}} from $\omega^\ell$, $(V_\ell,\bullet,\b1)$
is the {\em fusion algebra}.

Since the fusion algebra is well-known{\footnote{See, for example,
 \cite[Proposition 6.1]{beauville}.}}
to be
 semisimple, the CohFT with unit $\Omega^\ell$
 is also semisimple,  
The subject has a history starting in the mid 80s
 with the discovery and in 90s with several 
proofs of the {\em Verlinde formula} for the rank 
$d_g(\mu_1, \ldots, \mu_n)$, see \cite{beauville} for an overview.

Hence, steps (i) and (ii) of the computational strategy
of Section \ref{ww33} for $\Omega^\ell$ are complete (and have been
for many years). Step (iii) is the jump to moduli.

\subsection{Path to the $R$-matrix}\label{pprr}
The shifted $r$-spin CohFT $\widehat{\W}^r$ has an Euler field
obtained from the pure dimensionality of Witten's $r$-spin class which
was used to find the unique $R$-matrix in Section \ref{wittenr}. The CohFT $\Omega^\ell$
is not of pure dimension and has no Euler field.
A different path to the $R$-matrix is required here.

The restriction of the tensor $\Omega^\ell_{g,n}$ to the open set of
nonsingular curves
 $$\cM_{g,n}\subset \oM_{g,n}$$
forgets a lot of the  data of the CohFT. However, by \cite[Lemma 2.2]{MOPPZ},
the restriction is enough to uniquely determine the $R$-matrix
of $\Omega^\ell$. Fortunately, the restriction 
is calculable in closed form: 
\begin {itemize}
\item [$\bullet$] the first Chern class of the Verlinde bundle over 
$\cM_{g, n}$ 
is found in \cite {T}, 
\item [$\bullet$] the existence of a projectively flat connection{\footnote{Often
called the Hitchin connection.}} \cite{TUY}
on the Verlinde bundle over $\cM_{g, n}$  then determines
the full Chern character over $\mathcal M_{g, n}$.
\end {itemize}
As should be expected,
the computation of $\Omega^\ell$ relies significantly upon the past study
of the Verlinde bundles.

%Two basic invariants of the Verlinde bundles have been studied so far. First, $$\text{ rank } \mathbb E_{g} (\mu_1, \ldots, \mu_n)=d_g(\mu_1, \ldots, \mu_n)$$ is given by the Verlinde formula; see for instance \cite {beauville}. Second, an explicit closed expression for the first Chern class was obtained in genus $0$ in \cite{F},~\cite{Mu}, in genus $1$ and for one marking in~\cite{F}, and in arbitrary higher genus in~\cite {MOP}. We calculate here the total Chern character. As a consequence of our formulas, the higher Chern classes are seen to lie in the tautological ring
%$$RH^\star(\oM_{g,n}) \subset H^\star(\oM_{g,n})\, , $$ see~\cite{FP}.

\subsection {The $R$-matrix} For a simple Lie algebra $\mathfrak g$ and 
a level $\ell$, the conformal anomaly is 
 $$c=c(\mathfrak g, \ell)=\frac{\ell \dim \mathfrak g}{\check{h}+\ell}\, ,$$ where $\check{h}$ is the dual Coxeter number. For each  representation with highest weight $\mu$ of level $\ell$, define
$$
\mathsf w(\mu)=\frac{(\mu, \mu+2\rho)}{2(\check{h}+\ell)}.$$ Here, $\rho$ is half of the sum of the positive roots, and the Cartan-Killing form $(, )$ is 
normalized so that the longest root $\theta$ satisfies $$(\theta, \theta)=2.$$ 

\begin{example} For ${\mathfrak g} = {\mathfrak s}{\mathfrak l} (r, {\mathbb C}),$ the highest weight of a representation of level $\ell$ is given by an $r$-tuple of integers 
$$\mu = (\mu^{1}, \ldots, \mu^{r})\, , \, \, \ell \geq \mu^{1} \geq \cdots \geq \mu^{r} \geq 0\, ,$$ defined up to shifting the vector components by the same integer. 
Furthermore, we have
$$ c (\mathfrak g, \ell) = \frac{\ell (r^2-1)}{\ell+r},$$
$${\mathsf w} (\mu)=\frac{1}{2(\ell + r)}\left(\sum_{i=1}^{r} (\mu^i)^2-\frac{1}{r}\left(\sum_{i=1}^r \mu^i\right)^2+ \sum_{i=1}^{r} (r-2i+1)\mu^i\right).$$
\end{example}

Via the path to the $R$-matrix discussed in Section \ref{pprr}, a 
simple
closed form for the $R$-matrix of $\Omega^\ell$ is found in \cite{MOPPZ}
using 
the constants $c(\mathfrak g,\ell)$ and ${\mathsf{w}}(\mu)$ from
representation theory.

\begin{theorem}\label{moppz}
The CohFT $\Omega^\ell$ is reconstructed from the topological part
$\omega^\ell$ by the diagonal $R$-matrix
$$
{R(z)}_\mu^\mu = 
\exp \left(t z \cdot 
\left(-{\mathsf w}(\mu) + \frac{c({\mathfrak g}, \ell)}{24}\right) \right)\, .
$$
\end{theorem}

For the Lie algebra $\mathfrak g=\mathfrak{sl}_2$ at level $\ell=1$,
there are only two representations $\{ \emptyset, \square\}$
to consider{\footnote{Here, $\emptyset$
is the trivial representation (corresponding to $\b1$) and
$\square$ is the standard representation.}}, 
$$
c(\mathfrak{sl}_2,1) = 1, \quad \mathsf w(\emptyset)=0,\quad \mathsf w(\square)=\frac{1}{4}.
$$ 
As an example of the reconstruction result of Theorem \ref{moppz}, 
$$\Omega^\ell=R. \omega^\ell\, ,$$
the total Chern character
$\text{ch}\,\mathbb V_g(\square,\ldots,\square)$ at $t=1$ is
\begin{equation*}
\exp\left(-\frac{\lambda_1}{2}\right)\cdot 
\sum_{\Gamma\in \mathsf{G}_{g,n}^{\mathsf{even}}}
\frac{2^{g-h^1(\Gamma)}}{|\text{Aut}(\Gamma)|}\cdot 
\iota_{\Gamma\star} 
\left(\prod_{e\in \E} \frac{1-\exp\left(-\frac{1}{4}(\psi'_e+\psi''_e)\right)}{\psi_e'+\psi''_e}\cdot \prod_{l\in \L} e^{-\psi_l/4} \right)\, 
\end{equation*}
in $H^*(\oM_{g,n},\mathbb{Q})$. A few remarks about the above formula are required:
\begin {itemize} 
\item[$\bullet$] The classes $\lambda_1$ and $\psi$ are
the first Chern classes of the Hodge bundle and the
cotangent line bundle respectively.
\item [$\bullet$] The sum is over the set of {\em even} stable graphs, 
$$ \mathsf{G}_{g,n}^{\mathsf{even}} \subset \mathsf{G}_{g,n}\, ,$$
defined by requiring the valence $\mathsf{n}(v)$ to be even
for every vertex $v$ of the graph.
\item [$\bullet$] 
The Verlinde rank $d_{\mathsf{g}(v)}(\square, \ldots,\square)$ with 
$\mathsf{n}(v)$ insertions 
equals $2^g$ in the even case, see \cite{beauville}. 
The product of $2^{\mathsf{g}(v)}$ over the vertices of $\Gamma$ yields $2^{g-h^1(\Gamma)}$, where $h^1$ denotes the first Betti number.
\end{itemize} 

\subsection{Questions}

A different approach to the calculation of the 
Chern character of the Verlinde bundle in the
$\mathfrak{sl}_2$ case (for every level) was pursued in \cite{FMP}
using the geometry introduced by Thaddeus \cite{Thaddeus} to
prove the Verlinde formula. The outcome of \cite{FMP} is
a more difficult calculation (with a much more complex
answer), but with one advantage: the projective flatness of the 
Hitchin connection over $\cM_{g,n}$ is {\em not} used. When the flatness
is introduced, the method of \cite{FMP} yields tautological relations.
Unfortunately, no such relations are obtained by the
above $R$-matrix calculation of $\Omega^\ell$ since the 
projective flatness is an input.

\begin{question}\label{ff446} Is there an alternative
 computation of $\Omega^\ell$ which
does {\em not}
use the projective flatness of the Hitchin connection
and which systematically
produces tautological relations in $RH^*(\cM_{g,n},\mathbb{Q})?$ 
\end{question}

Of course, if the answer to Question \ref{ff446} is yes, then the next
question is whether {\em all} tautological relations are produced.{\footnote{I first heard an early version of this question from R. Bott at Harvard
 in the 90s.}}

\section{Gromov-Witten theory of $\Hilb$} \label{hhilb}

\subsection{$\T$-equivariant cohomology of $\Hilb$}
The {\em Hilbert scheme} $\Hilb$
of $m$ points in the plane $\CC^2$ 
parameterizes ideals $\cI\subset \CC[x,y]$ of colength $m$,
$$
\dim_\CC {\CC[x,y]}/{\cI} = m \,. 
$$
The Hilbert scheme $\Hilb$ is a nonsingular, irreducible,
quasi-projective variety 
of  
dimension $2m$,
see
\cite{Nak} for an introduction. 
An open dense set of $\Hilb$ parameterizes 
ideals associated to configurations of
$m$ distinct points.

The symmetries of $\CC^2$ lift to the Hilbert scheme. 
The algebraic torus 
$$\T=(\CC^*)^2$$ 
acts diagonally on $\CC^2$ by scaling coordinates,
$$
(z_1,z_2) \cdot (x,y) = (z_1 x, z_2 y)\, .
$$
We review the Fock space description of the $\T$-equivariant 
cohomology of the Hilbert scheme of points of $\CC^2$ following
the notation of \cite[Section 2.1]{op}. 

By definition, the {\em Fock space} $\cF$ 
is freely generated over $\Q$ by commuting 
creation operators $\alpha_{-k}$, $k\in\ZZ_{>0}$,
acting on the vacuum vector $\vac$. The annihilation 
operators $\alpha_{k}$, $k\in\ZZ_{>0}$, kill the vacuum 
$$
\alpha_k \cdot \vac =0,\quad k>0 \,,
$$
and satisfy the commutation relations
$\left[\alpha_k,\alpha_l\right] = k \, \delta_{k+l}$.

A natural basis of $\cF$ is given by 
the vectors  
\begin{equation}
  \label{basis}
  \lv \mu \rang = \frac{1}{\zz(\mu)} \, \prod_i \alpha_{-\mu_i} \, \vac \,
\end{equation}
indexed by partitions 
$\mu$. Here, $\zz(\mu)=|\Aut(\mu)| \, \prod_i \mu_i$ is the usual 
normalization factor. 
Let the length $\ell(\mu)$ denote the number of 
parts of the partition $\mu$.

The {\em Nakajima basis} defines a canonical isomorphism,
\begin{equation}
\cF \otimes _{\mathbb Q} {\mathbb Q}[t_1,t_2]\stackrel{\sim}{=} 
\bigoplus_{n\geq 0} H_{\T}^*(\Hilb,{\mathbb Q}).
\label{FockHilb}
\end{equation}
The Nakajima basis element corresponding to  
$\lv \mu \rang$  is
$$\frac{1}{\Pi_i \mu_i} [\mathcal{V}_\mu]$$
where $[\mathcal{V}_\mu]$ is (the cohomological dual of)
the class of the subvariety of $\mathsf{Hilb}^{|\mu|}(\CC^2)$
with generic element given by a union of 
schemes of lengths $$\mu_1, \ldots, \mu_{\ell(\mu)}$$ supported
at $\ell(\mu)$ distinct points{\footnote{The points and
parts of $\mu$ are considered
here to be unordered.}} of $\CC^2$. 
The vacuum vector $\vac$ corresponds to the unit in 
$$\mathsf{1}\in H_\T^*(\mathsf{Hilb}^0(\CC^2), {\mathbb Q})\, .$$
The variables $t_1$ and $t_2$ are the equivariant parameters corresponding
to the weights of the $\T$-action on the tangent
space 
$\text{Tan}_0(\CC^2)$ 
at the origin of $\CC^2$.

The subspace  $\cF_m\subset \cF\otimes_{\mathbb Q} {\mathbb Q}[t_1,t_2]$ corresponding to $H^*_\T(\Hilb,{\mathbb Q})$
is spanned by the vectors \eqref{basis} with $|\mu|=n$. 
The subspace can also 
be described as the $n$-eigenspace of the {\em  energy operator}: 
$$
|\cdot| = \sum_{k>0} \alpha_{-k} \, \alpha_k \,.
$$
The vector $\lv 1^n \rang$ corresponds to the unit 
$$\mathsf{1}\in H^*_\T(\Hilb,{\mathbb Q})\, .$$

%A straightforward calculation shows 
%\begin{equation}\label{DDDD}
%D = - \lv 2,1^{n-2} \rang \,.
%\end{equation}

The standard inner product on the $\T$-equivariant cohomology of
$\Hilb$ 
induces the following 
{\em nonstandard} inner product on Fock space after an extension of scalars:
\begin{equation}
  \label{inner_prod}
  \lang \mu | \nu \rang = 
\frac{(-1)^{|\mu|-\ell(\mu)}}{(t_1 t_2)^{\ell(\mu)}} 
\frac{\delta_{\mu\nu}}{\zz(\mu)} \,. 
\end{equation}
%
%With respect to the inner product, 
%%
%\begin{equation}
%  \label{adjoint}
%  \left(\alpha_{k}\right)^* = (-1)^{k-1} (t_1 t_2)^{\sgn(k)} \, 
%\alpha_{-k}\,.
%\end{equation}
%

\subsection{Gromov-Witten CohFT of $\Hilb$}

Let $m>0$ be a colength. Let $(V_m,\eta, \b1)$ be the
following triple:
\begin{enumerate}
\item[$\bullet$] $V_m$ is the free $\mathbb{Q}(t_1,t_2)[[q]]$-module 
$\cF_m\otimes_{\mathbb{Q}[t_1,t_2]}\mathbb{Q}(t_1,t_2)[[q]]$,
\item[$\bullet$] $\eta$ is the non-degenerate symmetric 2-form \eqref{inner_prod},

\item[$\bullet$] $\b1$ is the basis element $\lv 1^n\rang$. 
\end{enumerate}

Since $\Hilb$ is {\em not} proper, the Gromov-Witten
theory is only defined after localization by $\T$.
The CohFT with unit
$$\Omega^{\Hilb}=(\Omega_{g,n}^{\Hilb})_{2g-2+n>0}$$
is defined via the localized $\T$-equivariant Gromov-Witten classes of
$\Hilb$,
$$\Omega_{g,n}^{\Hilb} \in H^*(\oM_{g,n},\mathbb{Q}(t_1,t_2)[[q]]) \otimes 
(V_m^*)^n\,.$$
Here, $q$ is the Novikov parameter.
% with curve class degree
%defined by pairing with $D$,
%
Curves of degree $d$ are counted with 
weight $q^d$, where the curve degree is defined
by the pairing with the divisor $$D= -\lv 2,1^{m-2}\rang\,,\ \ \ \
 d=\int_\beta D\, .$$
%The ordinary multiplication in 
%$\T$-equivariant cohomology is recovered by
%setting $q=0$.

Formally, $\Omega^{\Hilb}$ is a CohFT {\em not} over the field $\mathbb{Q}$ as
in the $r$-spin and Verlinde cases, but
over the ring 
$\mathbb{Q}(t_1,t_2)[[q]]$.
%A basic result of \cite{OP,HHP}
%is {\em rationality} of the $q$-dependence: $\Omega^{\Hilb}$ is defined over
%the field $\mathbb{Q}(t_1,t_2,q)$. 
To simplify notation,
let $$\Omega^m= \Omega^{\Hilb}\, .$$

\subsection{Genus 0}
Since the $\T$-action on $\Hilb$ has finitely many $\T$-fixed
points, the localized $\T$-equivariant cohomology $$H^*_{\T}(\Hilb,\mathbb{Q})\otimes_{\mathbb{Q}[t_1,t_2]} \mathbb{Q}(t_1,t_2)$$
is semisimple. 
At $q=0$,
the quantum cohomology ring,
\begin{equation}
(V_m, \bullet, \b1)\, , \label{ffdd55}
\end{equation} 
defined by $\Omega^m$ specializes to the localized
$\T$-equivariant cohomology of $\Hilb$.
Hence, the quantum cohomology \eqref{ffdd55} is semisimple
over the ring $\Q(t_1,t_2)[[q]]$, see \cite{YPP}.

Let $\mathsf{M}_D$ denote the
operator of $\T$-equivariant quantum multiplication by the 
divisor $D$. 
A central result of \cite{op} is the following explicit formula for
 $\mathsf{M}_D$ an as operator on Fock space:
\begin{multline*}
%  \label{theM} 
\MM_D(q,t_1,t_2) = (t_1+t_2) \sum_{k>0} \frac{k}{2} \frac{(-q)^k+1}{(-q)^k-1} \,
 \alpha_{-k} \, \alpha_k\ \, -\, \frac{t_1+t_2}{2}\frac{(-q)+1}{(-q)-1}|\cdot |
 \\
 +  
\frac12 \sum_{k,l>0} 
\Big[t_1 t_2 \, \alpha_{k+l} \, \alpha_{-k} \, \alpha_{-l} -
 \alpha_{-k-l}\,  \alpha_{k} \, \alpha_{l} \Big] 
 \, .
\end{multline*}
The $q$-dependence of $\MM_D$ occurs only in the first two terms
(which act diagonally in the basis \eqref{basis}). 

Let $\mu^1$ and $\mu^2$ be partitions of $m$. The  
$\T$-equivariant Gromov-Witten invariants of  $\Hilb$ in genus $0$ 
with 3 cohomology insertions 
given (in the Nakajima basis) by $\mu^1$, $D$, and $\mu^2$ are determined by $\MM_D$:
$$\sum_{d=0}^\infty \Omega^m_{0,3,d} (\mu^1, D, \mu^2)\, q^d\ =\ 
\big\langle \mu^1\, \big| \, \MM_D  \, \big|\,  \mu^2 \big\rangle\, .$$
The following result is proven in \cite{op}.

\begin{theorem} \label{unnq}
The restriction of $\Omega^m$
to genus 0 is uniquely and effectively determined from the
calculation of $\MM_D$.
\end{theorem}

While Theorem \ref{unnq} in principle completes the genus 0 study of
$\Omega^m$, the result is not as strong as the  genus 0
determinations in the $r$-spin and Verlinde cases. The proof
of Theorem \ref{unnq} provides an effective linear algebraic
procedure, but not a formula,
for calculating the genus 0 part of $\Omega^m$ from $\MM_D$.

%By \cite[Corollary ?]{op},  
%The divisor class $D$ generates the small quantum ring
%$$QH^*_T(\Hilb,{\mathbb Q})$$ 
%over $\Q(q,t_1,t_2)$.  The small quantum ring structure
%is therefore determined. 

%Let $\mu^1,\ldots, \mu^r\in \text{Part}(n)$. 
%The $\T$-equivariant Gromov-Witten series in genus $g$,
%$$\blang \mu^1, \mu^2, \ldots, \mu^r \brang_{g}^{\Hilb}
%\, =\,  \sum_{d=0}^\infty 
%\blang \mu^1, \mu^2, \ldots, \mu^r \brang_{g,d}^{\Hilb} q^d
%\, \in \mathbb{Q}[[q]]\, ,$$
%is a sum over the degree $d$ with variable $q$.
%The $\T$-equivariant Gromov-Witten series in genus 0,
%$$\blang \mu^1, \mu^2, \ldots, \mu^r \brang_{0}^{\Hilb}
%= \sum_{d=0}^\infty 
%\blang \mu^1, \mu^2, \ldots, \mu^r \brang_{0,d}^{\Hilb} q^d\, ,
%$$
%can be calculated from the special 3-point
%invariants \eqref{v233}, see  \cite[Section 4.2]{op}. 

\subsection{The topological part}
Let $\omega^m$ be the topological part of the CohFT with unit $\Omega^m$.
A closed
formula for  $\omega^m$ can not be expected since 
closed formulas are already missing in the genus 0 study. 

The CohFT with unit $\omega^m$ has been considered earlier from another 
perspective. Using fundamental correspondences \cite{MNOP},
$\omega^m$ is equivalent to the {\em local GW/DT theory} of 3-folds of
the form
\begin{equation}\label{gg3311}
\C^2 \times C\, ,
\end{equation}
where $C$ is a curve or arbitrary genus. Such local
theories have been studied extensively \cite{bp}  in the
investigation of the GW/DT theory of 3-folds{\footnote{
A natural generalization of the geometry \eqref{gg3311} is
to consider the 3-fold total space, 
$$L_1 \oplus L_2 \rightarrow C\, ,$$
of a sum of  line bundles $L_1,L_2 \rightarrow C$. For
particular pairs $L_1$ and $L_2$, simple closed form solutions
were found \cite{bp}
and have later played a role in the study of the structure
of the Gromov-Witten theory of Calabi-Yau 3-folds by Ionel and Parker \cite{IP}.}}.

\subsection{The $R$-matrix}
Since $\Omega^m$ is not of pure dimension (and does not carry an
Euler field),
the $\mathsf{R}$-matrix is {\em not} determined by the
$\T$-equivariant genus $0$ theory alone.
As in the Verlinde case, a different method is required. 
Fortunately,  together with the divisor equation, an evaluation of the
$\T$-equivariant higher genus theory in degree $0$ is enough
to uniquely determine the $\mathsf{R}$-matrix. 

Let
$\text{Part}(m)$ be the set of of partitions
of $m$ corresponding to the $\T$-fixed
points of $\Hilb$. For each $\eta\in \text{Part}(m)$,  let
 $\text{Tan}_\eta(\Hilb)$ be the $\T$-representation on
 the tangent space at the $\T$-fixed point
 corresponding to $\eta$.
As before, let
$$\mathbb{E}\rightarrow \oM_{g,n}$$
be the Hodge bundle. The follow result is proven in \cite{HHP}.

\begin{theorem}\label{ff11} The $\mathsf{R}$-matrix 
of $\Omega^m$
 is uniquely determined by the divisor equation
and the degree 0 invariants
\begin{eqnarray*}
\big\langle \mu \big\rangle_{1,0}^{\Hilb} &=& \sum_{\eta\in \text{\em Part}(m)} \mu|_\eta \int_{\oM_{1,1}}
\frac{e\left(\mathbb{E}^* \otimes \text{\em Tan}_\eta(\Hilb)\right)}
{e\left(\text{\em Tan}_\eta(\Hilb)\right)
}\, , 
\\
\big\langle \, \big\rangle_{g\geq 2,0}^{\Hilb} &=& \sum_{\eta\in \text{\em Part}(m)} \int_{\oM_g}
\frac{e\left(\mathbb{E}^* \otimes \text{\em Tan}_\eta(\Hilb)\right)}
{e\left(\text{\em Tan}_\eta(\Hilb)\right)
}\,
.
\end{eqnarray*}
\end{theorem}

While 
Theorem \ref{ff11} is weaker than the explicit $R$-matrix
solutions in the $r$-spin and Verlinde cases, the
result nevertheless has several consequences. The first 
is a rationality result \cite{HHP}.

\begin{theorem}\label{crcr} For all genera $g\geq 0$ and 
$\mu^1,\ldots,\mu^n\in \text{\em Part}(m)$, 
the series{\footnote{As always,
$g$ and $n$ 
are required
to be in the stable range $2g-2+n>0$.}}
$$\int_{\oM_{g,n}}\Omega^m_{g,n}( \mu^1, \ldots, \mu^n)\, \in\, \mathbb{Q}(t_1,t_2)[[q]]$$
is the Taylor expansion in $q$ of a rational
function in $\mathbb{Q}(t_1,t_1,q)$.
\end{theorem}

The statement of Theorem \ref{crcr} can be strengthened (with 
an $R$-matrix argument using Theorem \ref{ff11}) to prove
that  the CohFT with unit $\Omega^m$ can be defined over the field $\mathbb{Q}(t_1,t_2,q)$.

\subsection{Crepant resolution}
The Hilbert scheme of points of $\mathbb{C}^2$  is well-known to be a crepant resolution
of the symmetric product,
$$\epsilon:\Hilb \ \rightarrow\ {\Sym}(\mathbb{C}^2)= (\mathbb{C}^2)^m/S_m\, .$$
Viewed as an {\em orbifold}, the symmetric product $\Sym(\mathbb{C}^2)$ has
a $\T$-equivariant
Gromov-Witten theory with insertions indexed by partitions of $m$
and an associated CohFT with unit $\Omega^{{\Sym}(\mathbb{C}^2)}$ 
determined by the Gromov-Witten classes.
The CohFT $\Omega^{{\Sym}(\mathbb{C}^2)}$ is defined over the ring
$\mathbb{Q}(t_1,t_2)[[u]]$, where $u$ is variable associated to
the free ramification points, see \cite{HHP} for a detailed treatment.

In genus $0$, the equivalence of the $\T$-equivariant
Gromov-Witten theories of $\Hilb$ and the orbifold ${\Sym}(\mathbb{C}^2)$
was proven{\footnote{The
prefactor
$(-i)^{\sum_{i=1}^n \ell(\mu^i)-|\mu^i|}$
was treated incorrectly in
\cite{bg} because of an
arthimetical error. The
prefactor here is correct.}}
in \cite{bg}. Another consequence of the $R$-matrix study of
$\Omega^{\Hilb}$ is the proof in \cite{HHP}
of the {\em crepant resolution conjecture} here.

\begin{theorem}\label{crc} For all genera $g\geq 0$ and 
$\mu^1,\ldots,\mu^n\in \text{\em Part}(m)$, we have  
$$\Omega_{g,n}^{\Hilb}(\mu^1, \ldots, \mu^n) =
(-i)^{\sum_{i=1}^n \ell(\mu^i)-|\mu^i|}\, 
\Omega^{{\Sym}(\mathbb{C}^2)}_{g,n}(\mu^1, \ldots, \mu^n)$$
 after the variable change $-q=e^{iu}$.
\end{theorem}

The variable change of Theorem \ref{crc} is well-defined by the rationality
of Theorem \ref{crcr}. The analysis \cite{opQ}
of the quantum differential equation
of $\Hilb$ plays an important role in the proof. 
Theorem \ref{crc} is closely related
to the GW/DT correspondence for local curves \eqref{gg3311} 
in families, see \cite{HHP}.

\subsection{Questions}

The most basic open question is to find an expression for
the $R$-matrix of $\Omega^m$ in terms of natural operators on
Fock space.

\begin{question}\label{nsns} Is there a representation theoretic formula for the
$R$-matrix of the CohFT with unit $\Omega^m$?
\end{question}

The difficulty in attacking Question \ref{nsns}
starts with the lack of higher genus 
calculations in closed form.
The first
nontrivial example \cite{HHP} occurs  in genus 1
for the Hilbert scheme of 2 points:
\begin{equation}\label{dxx12}
\int_{\oM_{1,1}} \Omega^{\mathsf{Hilb}^2(\mathbb{C}^2)}_{1,1}\big((2)\big)= -\frac{1}{24}\frac{(t_1+t_2)^2}{t_1t_2} \cdot \frac{1+q}{1-q}\, .
\end{equation}
While there are numerous calculations to do, the higher $m$
analogue of \eqref{dxx12} surely has a simple answer.

\begin{question}
Calculate the series 
\begin{equation*}
\int_{\oM_{1,1}} \Omega^{\mathsf{Hilb}^m(\mathbb{C}^2)}_{1,1}\big((2,1^{n-2})\big)\
\in\ \mathbb{Q}(t_1,t_2,q)\, ,
\end{equation*}
in closed form for all $m$.
\end{question}

\vspace{+16 pt}
\noindent Departement Mathematik \\
\noindent ETH Z\"urich \\
\noindent rahul@math.ethz.ch

\end{document}